\documentclass[12pt]{article}

\usepackage[margin=1.5cm]{geometry}

\usepackage{cite}

\usepackage{tikz}

\usepackage{bm}

\usepackage{mathtools}
\usepackage{graphicx}
\usepackage{amssymb,amsmath,wasysym}
\usepackage{amsfonts}
\usepackage{comment}
\usepackage{enumerate}

\newtheorem{Example}{Example}
\newtheorem{Remark}{Remark}

\newcommand{\bpi}{\mbox{\boldmath$\pi$}}

\newcommand{\balpha}{\mbox{\boldmath$\alpha$}}

\newcommand{\btheta}{\mbox{\boldmath$\theta$}}
\newcommand{\bepsilon}{\mbox{\boldmath$\epsilon$}}

\begin{document}
	

\title{Sensitivity analysis of Quasi-Birth-and-Death processes}

\author{
%
%
Anna Aksamit
\thanks{Anna Aksamit is supported by the Australian Research Council Early Career Researcher Award DE200100896.}
~~\thanks{School of Mathematics and Statistics, The University of Sydney, NSW 2006, Australia, email: anna.aksamit@sydney.edu.au}
\and 
Ma{\l}gorzata M. O'Reilly
\thanks{Ma{\l}gorzata M. O'Reilly is supported by the Australian Research Council Discovery Project DP180100352.}
~~\thanks{Discipline of Mathematics, University of Tasmania, Tas 7001, Australia, email: malgorzata.oreilly@utas.edu.au}
\and
Zbigniew Palmowski
\thanks{Zbigniew Palmowski was partially supported by the National Science Centre (Poland) under the grant 2021/41/B/HS4/00599.}
~~\thanks{Faculty of Pure and Applied Mathematics, Wroc{\l}aw University of Science and Technology, 50-370 Wroc{\l}aw, Poland, email: zbigniew.palmowski@pwr.edu.pl}
}

\date{\today}

\maketitle

\section{Introduction}

Quasi-birth-and-death processes (QBDs) is the fundamental class of Markovian models in the theory of matrix-analytic methods, with a level-variable $X(t)$ and a phase-variable $\varphi(t)$ forming a two-dimensional state space. In many applications of the QBDs, the phase variable $\varphi(t)$ is used to model information about the underlying environment that drives the evolution of some system.

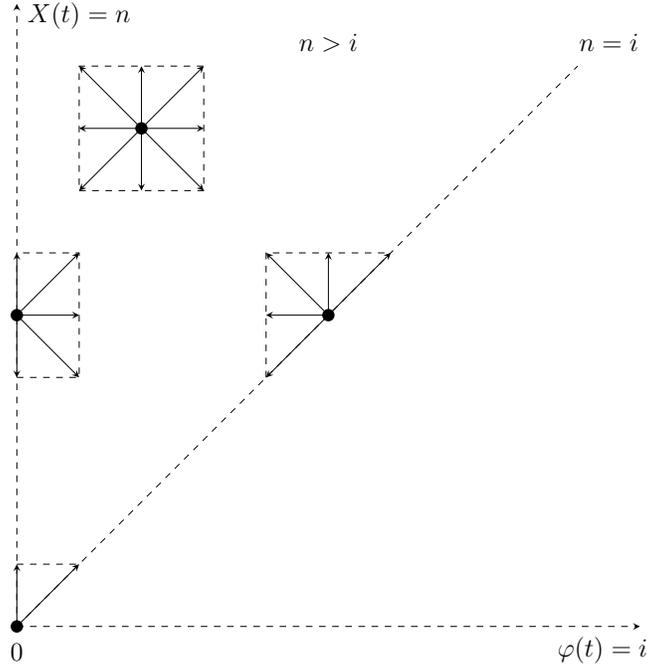
\begin{figure}[h!]
	\begin{center}
		\resizebox{0.48\textwidth}{!}{
			\begin{tikzpicture}[>=stealth,redarr/.style={->}]
			

%
%
%
%
%
%
			
			\draw [ dashed, ->] (0,0) -- (10,0);
			\draw [ dashed, ->] (0,0) -- (0,10);
			\draw [dashed] (0,0) -- (9,9);
			
			\draw (1,10) node[anchor=north, below=-0.17cm] {{\color{black} $X(t)=n$}};			
			
			\draw (9.4,-0.2) node[anchor=north, below=-0.17cm] {{\color{black} $\varphi(t)=i$}};			
			
			\draw (0,-0.3) node[anchor=north, below=-0.17cm] {{\color{black} $0$}};
			
			\draw (9.5,9.5) node[anchor=north, below=-0.17cm] {{\color{black} $n=i$}};
			
			\draw (5,9.5) node[anchor=north, below=-0.17cm] {{\color{black} $n>i$}};			
			


			\node at (2,8) [black,circle,fill,inner sep=2pt]{};
			\draw [black, ->] (2,8) -- (3,8);
			\draw [black, ->] (2,8) -- (1,8);
			\draw [black, ->] (2,8) -- (2,9);
			\draw [black, ->] (2,8) -- (2,7);
			\draw [black, ->] (2,8) -- (3,9);
			\draw [black, ->] (2,8) -- (1,7);
			\draw [black, ->] (2,8) -- (3,7);
			\draw [black, ->] (2,8) -- (1,9);
			\draw [dashed] (1,7) -- (3,7) -- (3,9) -- (1,9) -- (1,7);

			\node at (5,5) [black,circle,fill,inner sep=2pt]{};
			\draw [black, ->] (5,5) -- (4,5);
			\draw [black, ->] (5,5) -- (5,6);
			\draw [black, ->] (5,5) -- (6,6);
			\draw [black, ->] (5,5) -- (4,4);
			\draw [black, ->] (5,5) -- (4,6);
			\draw [dashed] (4,4) -- (4,6) -- (6,6);


			\node at (0,5) [black,circle,fill,inner sep=2pt]{};
			\draw [black, ->] (0,5) -- (0,6);
			\draw [black, ->] (0,5) -- (0,4);
			\draw [black, ->] (0,5) -- (1,6);
			\draw [black, ->] (0,5) -- (1,4);
			\draw [black, ->] (0,5) -- (1,5);
			\draw [dashed] (0,4) -- (1,4) -- (1,6) -- (0,6);
			
			\node at (0,0) [black,circle,fill,inner sep=2pt]{};
			\draw [black, ->] (0,0) -- (0,1);
			\draw [black, ->] (0,0) -- (1,1);
			\draw [dashed] (0,1) -- (1,1);

			\end{tikzpicture}
		}
	\end{center}
	\caption[Evolution of a QBD]{Evolution of a QBD $\{(\varphi(t),X(t)):t\geq 0\}$ modelling the total number of patients $X(t)$ and some information $\varphi(t)$ about the system e.g. (here) the number of patients of a particular class such as complex patients, whose treatment is likely to require more resources and time~\cite{Gus2021}.}
	\label{fig:RWquadrant}
\end{figure}

A QBD is a model that lends itself to representing healthcare system in a natural, intuitive manner (see Figure~\ref{fig:RWquadrant}), and so the application potential of the QBDs in this area is immense, as demonstrated by Heydar et al.~\cite{2021HOTFTT} and Grant~\cite{Gus2021}. As another example of application in real world systems, QBDs have been applied in the analysis of evolution of gene families e.g. in Diao et al.~\cite{DSLOH2020}.
 

An initial sensitivity analysis of a level-dependent QBD (LD-QBD) has been performed by G{\'o}mez-Corral and L{\'o}pez-Garc{\'i}a in~\cite{2018CL}. Here, we build on these ideas and extend the analysis to a wide range of metrics of interest, with a particular focus on applications in modelling healthcare systems.

Suppose that the generator ${\bf Q}(\btheta)$ of a LD-QBD depends on some parameters recorded in a row vector $\btheta=[\theta_i]_{i=1,\ldots,k}$. Various stationary (long-run) and transient (time-dependent) quantities in the analysis of the LD-QBDs, recorded as matrices ${\bf A}(\btheta)=[A_{ij}(\btheta)]$, can then be expressed using expressions involving ${\bf Q}(\btheta)$, and so, also depend on $\btheta$.

We derive the sensitivity analysis of relevant quantities, where given a vector of parameters $\btheta=[\theta_i]_{i=1,\ldots,k}$ and matrix ${\bf A}(\btheta)$, we write
\begin{eqnarray*}
	\frac{\partial {\bf A}(\btheta)}{\partial \btheta}
	&=&
	\left[
	\begin{array}{ccc}
		\frac{\partial {\bf A}(\btheta)}{\partial \theta_1}
		& \ldots
		& \frac{\partial {\bf A}(\btheta)}{\partial \theta_k}
	\end{array}
	\right].
\end{eqnarray*}
Here, we focus on the development of the key building blocks of the methodology. Full details of this work, including numerical examples, will be presented in our future paper.

\section{LD-QBD model}

Suppose that $\{(X(t),\varphi(t)):t\geq 0\}$ is a continuous-time Markov chain with a two-dimensional state $(X(t),\varphi(t))$ consisting of the level variable $X(t)$ and the phase variable, $\varphi(t)$, taking values in an irreducible state space given by
\begin{equation*}
\mathcal{S}=\{(n,i): n=0,1,2,\ldots,N ; i=0,1,\ldots,m_n\},
\end{equation*}
and with transition rates recorded in the generator matrix ${\bf Q}=[q_{(n,i)(n',i')}]_{(n,k),(n',i')\in\mathcal{S}}$ made of block matrices ${\bf Q}^{[n,n']}=[q_{(n,i)(n',i')}]_{i=0,1,\ldots,m_n,i'=0,1,\ldots,m_{n'})}$ such that
\begin{eqnarray*}
	\lefteqn{
		{\bf Q}
		=
		[{\bf Q}^{[n,n']}]_{n,n'=0,1,\ldots,N}
	}
	\\[1ex]
	&=&
	\begin{bmatrix}
		{\bf Q}^{[0,0]} & {\bf Q}^{[0,1]} & {\bf 0} & \cdots &  \cdots & {\bf 0}\\
		{\bf Q}^{[1,0]} & {\bf Q}^{[1,1]} & {\bf Q}^{[1,2]} & \cdots & \cdots & {\bf 0}\\
		\vdots & \vdots & \vdots & \cdots  & \cdots & \vdots\\
		{\bf 0} & {\bf 0} & {\bf 0} & \cdots & {\bf Q}^{[N,N-1]} & {\bf Q}^{[N,N]}
	\end{bmatrix},\nonumber\\
\end{eqnarray*}
so that only transitions to the neighbouring levels are possible. We refer to such process as a continuous-time level-dependent quasi-birth-and-death process (LD-QBD), that is bounded from above by level $N$. The level variable $X(t)$ may be used to record the number of individuals in some system at time $t$, while $\varphi(t)$ may be used to record some additional information about the system at time $t$.

Since a LD-QBD is a continuous-time Markov chain, standard expressions from the theory of Markov chains apply. However, as its state space $\mathcal{S}$ may be very large, here we apply ideas from the theory of matrix-analytic methods, which leads to efficient computational methods. We also consider the LD-QBD which does not have an upper boundary $N$.

We note that a range of various transient and stationary performance measures of such defined LD-QBD can be readily derived using the existing methods in the literature of matrix-analytic methods, see Ramaswami~\cite{ram1997}, Joyner and Fralix~\cite{joyner2016new}, and Phung-Duc et al.~\cite{phung2010simple}. Since these performance measures depend on the parameters of the model, we are interested in the sensitivity analysis of these measures.

%
%

\section{Quantities of interest}

We consider the following key quantities in the analysis of the LD-QBDs,
\begin{itemize}
	\item the long-run proportion of times spent in states $(n,i)$; 
	\item the distribution of times spent within levels contained in some set $\mathcal{A}\subset\{0,1,\ldots,N\}$ in a sample path for the process to first reach level $n\pm k$ and do so in state $(n\pm k,j)$ given start from state $(n,i)$;
	\item the distribution of the process observed at time $t$ given start from state $(n_0,i)$.
\end{itemize}
We follow the approach summarised in Grant~\cite{Gus2021}, which is built on the results in Ramaswami~\cite{ram1997}, Joyner and Fralix~\cite{joyner2016new}, and Phung-Duc et al.~\cite{phung2010simple}. We also state expressions for the relevant Laplace-Stieltjes transforms (LSTs) of various quantities. These can be inverted it using numerical inversion techniques by Abate and Whitt~\cite{abate1995numerical}, Den Iseger~\cite{DenIseger_2006}, or Horv{\'a}th et al.~\cite{horvath2020numerical}, to compute the corresponding quantities.

\subsection{Stationary distribution}\label{sec:StatDistr}

For all $n=0,1,\ldots,N$, $i=0,1,\ldots,m_n$, define the limiting probabilities
\begin{eqnarray*}
	\pi_{(n,i)}=\lim_{t\to\infty}\mathbb{P}\left(X(t)=n,\varphi(t)=i\right)
\end{eqnarray*}
recording the long-run proportions of time spent in states $(n,i)$, and collect these in a row vector $\bm{\pi}=[\bm{\pi}_n]_{n=0,1,\ldots,N}$, where $\bm{\pi}_n=[\pi_{(n,i)}]_{i=0,1,\dots,I_n}$.

To evaluate $\bm{\pi}$, we follow the approach in Grant~\cite{Gus2021}, and write the expressions for $\bm{\pi}_n$ in terms of $\bm{\pi}_N$ (rather than in terms of $\bm{\pi}_0$, since potential close-to-zero values in $\bm{\pi}_0$ may lead to computational errors). We consider matrices $$\widehat{\bf R}^{(n)}=[\widehat R^{(n)}_{ij}]_{i=0,1,\dots,m_{n-1},j=0,1,\dots,m_{n}}$$ recording  the expected times $\widehat R^{(n)}_{ij}$ spent in states $(n,j)$ per unit time spent in $(n+1,i)$, before returning to level $n+1$, given the process starts in state $(n+1,i)$. 

For $n=0,1,\ldots,N-1$, we apply the recursion,
\begin{eqnarray*}
	\widehat {\bf R}^{(0)}(s) &=& -{\bf Q}^{[1,0]}({\bf Q}^{[0,0]}-s{\bf I})^{-1},\\
	\widehat {\bf R}^{(n)}(s) &=& -{\bf Q}^{[n+1,n]}(\widehat {\bf R}^{(n-1)}(s){\bf Q}^{[n-1,n]}+{\bf Q}^{[n,n]}-s{\bf I})^{-1},  
\end{eqnarray*}
and then, with $\widehat{\bf R}^{(n)}=\widehat{\bf R}^{(n)}(0)$, let
\begin{eqnarray*}
\bm{\pi}_{n}&=&\bm{\pi}_{n+1} \widehat{\bf R}^{(n)}
=
\bm{\pi}_N \prod_{k=N-1}^{n} \widehat{\bf R}^{(k)},
\label{eq:pinRn}
\end{eqnarray*}
where $\bm{\pi}_N$ is the solution of the set of equations,
\begin{eqnarray}
\bm{\pi}_N
\left(
\widehat {\bf R}^{(N-1)}
{\bf Q}^{[N-1,N]}+{\bf Q}^{[N,N]} 
\right)\
&=& {\bf 0},
\nonumber\\
\bm{\pi}_N
\left(
{\bf 1}
+
\sum_{n=0}^{N-1}
\prod_{k=N-1}^{n}\widehat{\bf R}^{(k)}{\bf 1}
\right)&=&1.
\nonumber
\label{eq:piN}
\end{eqnarray}

\begin{Remark}
	Similar methods may be applied to derive $\bpi$ for a LD-QBD in which the level variable $X(t)\geq 0$ has no upper boundary $N$. First, following Phung-Duc et al.~\cite{phung2010simple}, find a sufficiently large truncation level $L$ such that the condition $||\widehat {\bf R}_{L}^{(n)}-\widehat {\bf R}_{L-1}^{(n)} ||<\epsilon$ is met for a required criterion $\epsilon>0$, where $\widehat {\bf R}_{L}^{(n)}$ and $\widehat {\bf R}_{L-1}^{(n)}$ denote matrix $\widehat {\bf R}^{(n)}$ computed for a LD-QBD with an upper boundary $N=L$ and $N=L-1$, respectively. Next, apply the approximation $\widehat {\bf R}^{(n)}\approx \widehat {\bf R}_{L}^{(n)}$. This technique may be applied for the remaining quantities.
\end{Remark}

\subsection{Sojourn times in specified sets}\label{sec:Sojtimes}

Let $\mathcal{A}\subset\{0,1,\ldots,N\}$ be some set of desirable or undesirable levels, such as $\mathcal{A}=\{0,\ldots,A\}$ or $\mathcal{A}=\{B,\ldots,N\}$ where $A$ and $B$ are some desirable or undesirable thresholds.

Let $\theta_n=\inf \{t>0:X(t)=n\}$ and $I(\cdot)$ be an indicator function. For any $n'\not=n$, let $${\bf W}^{n,n'}(s)=[W_{ij}^{n,n'}(s)]_{i=1,\ldots,m_n;j=1,\ldots,m_{n'}}$$ be a matrix such that the entry 
\begin{eqnarray*}
W_{ij}^{n,n'}(s)&=&
\mathbb{E}(
e^{-s\theta_{n'}}\times I\left(\varphi(\theta_{n'}\leq t,\theta_{n'}=j\right)
\\
&&\quad
 \ | \ X(0)=n,\varphi(0)=i))
\end{eqnarray*}
is the Laplace-Stieltjes transform (LST) of the time for the process to first visit level $n'$ and do so in phase $j$, given start from level $n$ in phase $i$.

Let ${\bf W}_{\mathcal{A}}^{n,n'}(s)=[W_{\mathcal{A};ij}^{n,n'}(s)]_{i=1,\ldots,m_n;j=1,\ldots,m_{n'}}$ be the LST matrix of the total time spent in the set $\mathcal{A}$ during a sample path corresponding to $W_{ij}^{n,n'}$.

Denote ${\bf G}^{n,n'}(s)={\bf W}^{n,n'}(s)$ and ${\bf G}_{\mathcal{A}}^{n,n'}(s)={\bf W}^{n,n'}(s)$ whenever $n'<n$, and ${\bf H}^{n,n'}(s)={\bf W}^{n,n'}(s)$ and ${\bf H}_{\mathcal{A}}^{n,n'}(s)={\bf W}_{\mathcal{A}}^{n,n'}(s)$ whenever $n'>n$.

We note that when $\mathcal{A}=\{n'+1,\ldots,N\}$, then clearly ${\bf G}_{\mathcal{A}}^{n,n'}(s)={\bf G}^{n,n'}(s)$. When $\mathcal{A}=\{0,\ldots,n'-1\}$, then ${\bf H}_{\mathcal{A}}^{n,n'}(s)={\bf H}^{n,n'}(s)$.

By standard decomposition of a sample path~\cite{Gus2021,ram1997,joyner2016new,phung2010simple,SOB2020}, we have,
\begin{eqnarray*}
{\bf G}_{\mathcal{A}}^{n,n-k}(s)&=&
{\bf G}_{\mathcal{A}}^{n,n-1}(s)
\times
\cdots
\times
{\bf G}_{\mathcal{A}}^{n-k+1,n-k}(s),
\end{eqnarray*}
where
\begin{eqnarray*}
{\bf G}_{\mathcal{A}}^{N,N-1}(s)&=&
-({\bf Q}^{[N,N]}-s{\bf I}\times I(N\in\mathcal{A}))^{-1}
{\bf Q}^{[N,N-1]},
\end{eqnarray*} 
and for $n=N-1,\ldots,1$,
\begin{eqnarray*}
{\bf G}_{\mathcal{A}}^{n,n-1}(s)
&=&
-\Big({\bf Q}^{[n,n]}-s{\bf I}\times I(n\in\mathcal{A})
\\
&&
\quad \quad 
+{\bf Q}^{[n,n+1]}{\bf G}_{\mathcal{A}}^{n+1,n}(s)
\Big)^{-1}
{\bf Q}^{[n,n-1]}.
\end{eqnarray*}

Similarly,
\begin{eqnarray*}
{\bf H}_{\mathcal{A}}^{n,n+k}(s)&=&
{\bf H}_{\mathcal{A}}^{n+1,n+2}(s)
\times
\cdots
\times
{\bf H}_{\mathcal{A}}^{n+k-1,n+k}(s),
\end{eqnarray*}
where
\begin{eqnarray*}
{\bf H}_{\mathcal{A}}^{0,1}(s)&=&
-({\bf Q}^{[0,0]}-s{\bf I}\times I(0\in\mathcal{A}))^{-1}
{\bf Q}^{[0,1]},
\end{eqnarray*} 
and for $n=N-1,\ldots,1$,
\begin{eqnarray*}
{\bf H}_{\mathcal{A}}^{n,n+1}(s)&=&
-\Big({\bf Q}^{[n,n]}-s{\bf I}\times I(n\in\mathcal{A})
\\
&&
\quad\quad 
+{\bf Q}^{[n,n-1]}{\bf H}^{n-1,n}(s)
\Big)^{-1}
{\bf Q}^{[n,n+1]}.
\end{eqnarray*}

\subsection{Distribution at time $t$}\label{sec:Distrt}

Suppose that the QBD starts from some level $n_0$ in some phase $i_0=0,1,\ldots,m_{n_0}$ according to the initial distribution of phases $\bm{\alpha}=[\alpha_j]_{j\in 0,1,\ldots,m_{n_0}}$ such that $\alpha_j=\mathbb{P}(\varphi(0)=j)$.

Define vector ${\bf f}(t)=[{\bf f}_n(t)]_{n=0,1,\ldots,N}$ such that
\begin{eqnarray*}
[{\bf f}_n(t)]_j&=& \mathbb{P}_{\bm{\alpha}}(X_t=n,\varphi(t)=j)
\end{eqnarray*}
is the probability that at time $t$ the process is on level $n$ and in phase $j$, given $\bm{\alpha}$; and the corresponding Laplace Transform vector $\widetilde{\bf f}(s)=[\widetilde{\bf f}_n(s)]_{n=0,1,\ldots,N}$ such that
\begin{eqnarray*}
[\widetilde{\bf f}_n(s)]_j&=& \int_{t=0}^{\infty} e^{-st} \mathbb{P}_{\bm{\alpha}}(X_t=n,\varphi(t)=j)dt
.
\end{eqnarray*} 

The Kolmogorov differential equations of the process are,
\begin{eqnarray*}
\frac{\partial {\bf f}_0(t)}{ \partial t}
&=&
{\bf f}_0(t){\bf Q}^{[0,0]}
+
{\bf f}_1(t){\bf Q}^{[1,0]}
,
\nonumber\\
\frac{\partial {\bf f}_n(t)}{ \partial t}
&=&
{\bf f}_n(t){\bf Q}^{[n,n]}
+
{\bf f}_{n-1}(t){\bf Q}^{[n-1,n]}
+
{\bf f}_{n+1}(t){\bf Q}^{[n+1,n]},
\nonumber\\
&&
\quad \mbox{ for }n=1,\ldots,N-1,
\nonumber\\
\frac{\partial {\bf f}_N(t)}{ \partial t}
&=&
{\bf f}_N(t){\bf Q}^{[N,N]}
+
{\bf f}_{N-1}(t){\bf Q}^{[N-1,N]},
\end{eqnarray*}
with the initial condition ${\bf f}_{n}(0)=\balpha I(n= n_0)$.

To evaluate $\widetilde{\bf f}(s)$, we apply the following recursion summarised in Grant~\cite{Gus2021},
\begin{eqnarray*}
\widetilde{\bf f}_{n_0}(s)&=& 
-\bm{\alpha}
\Big(
({\bf Q}^{[n_0,n_0]}-s{\bf I})
+
\widehat {\bf R}^{(n_0-1)}(s){\bf Q}^{[n_0-1,n_0]}
\\
&&
\quad\quad
+
\widetilde{\bf R}^{(n_0+1)}(s){\bf Q}^{[n_0+1,n_0]}
\Big)^{-1}
,
\end{eqnarray*}
and then for $n>n_0$,
\begin{eqnarray*}
\widetilde{\bf f}_{n}(s)&=& 
-\bm{\alpha}
{\bf H}^{n_0,n}(s)
\Big(
({\bf Q}^{[n,n]}-s{\bf I})
+
\widehat {\bf R}^{(n-1)}(s){\bf Q}^{[n-1,n]}
\\
&&
\quad\quad
+
\widetilde{\bf R}^{(n+1)}(s){\bf Q}^{[n+1,n]}
\Big)^{-1}
,
\end{eqnarray*}
and for $n<n_0$,
\begin{eqnarray*}
\widetilde{\bf f}_{n}(s)&=& 
-\bm{\alpha}
{\bf G}^{n_0,n}(s)
\Big(
({\bf Q}^{[n,n]}-s{\bf I})
+
\widehat {\bf R}^{(n-1)}(s){\bf Q}^{[n-1,n]}
\\
&&
+
\widetilde{\bf R}^{(n+1)}(s){\bf Q}^{[n+1,n]}
\Big)^{-1}
,
\end{eqnarray*}
where, for $n=N,\ldots,1$, we apply the recursion,
\begin{eqnarray*}
	\widetilde R_{N}({\bf 0}) &=& -{\bf Q}^{[N-1,N]}({\bf Q}^{[N,N]}-s{\bf I})^{-1},\\ 
\widetilde{\bf R}^{(n)}(s) &=& -{\bf Q}^{[n-1,n]}({\bf Q}^{[n,n]}-s{\bf I}+\widetilde{\bf R}^{(n+1)}(s){\bf Q}^{[n+1,n]})^{-1}.  
\end{eqnarray*}

\section{Sensitivity analysis}


First, we present simple examples to motivate the theory. 
\begin{Example}
	Suppose that $X(t)$ records the total number of customers in the system at time $t$. Let $\varphi(t)\in\{1,\ldots,k\}$ be the phase of the environment that drives the evolution of the system, so that $\lambda_{\varphi(t)}>0$ is the arrival rate to system (provided $X(t)<N$), and $\mu_{\varphi(t)}>0$ is the service rate per customer (provided $X(t)>0$). Assume that $\{\varphi(t):t\geq 0\}$ is a continuous-time Markov chain with generator ${\bf T}=[T_{ij}]$. The system can be modelled as a LD-QBD with generator ${\bf Q}(\btheta)=[q(\btheta)_{(n,i)(m,j)}]$ that depends on the vector of parameters $\btheta=[\lambda_1,\ldots,\lambda_k,\mu_1,\ldots,\mu_k]$, such that the nonzero off-diagonals $q(\btheta)_{(n,i)(m,j)}$ are given by
	\begin{eqnarray*}
		\left\{
		\begin{array}{ll}
			T_{ij}& j\not= i, m=n;\\
			\lambda_i& j = i, m=n+1, n<N;\\
			n\mu_i& j = i, m=n-1, n>0.
		\end{array}
		\right.
	\end{eqnarray*}
	
	Then $\frac{\partial}{\partial\btheta}{\bf Q}(\btheta)$ is given by
	\begin{eqnarray*}
		\frac{\partial q(\btheta)_{(n,i)(m,j)}}{\partial\lambda_i}&=&
		\left\{
		\begin{array}{ll}
			1& j = i, m=n+1, n<N;\\
			-1& j = i, m=n, 0\leq n<N;\\
			0& \mbox{otherwise;}	
		\end{array}
		\right.
	\end{eqnarray*}
	and
	\begin{eqnarray*}
		\frac{\partial q(\btheta)_{(n,i)(m,j)}}{\partial\mu_i}&=&
		\left\{
		\begin{array}{ll}
			n& j = i, m=n-1, n>0;\\
			-n& j = i, m=n, 0<n\leq N;\\
			0& \mbox{otherwise.}	
		\end{array}
		\right.
	\end{eqnarray*}
\end{Example}

\begin{Example}
	Suppose that customers of type $k\in\{1,2\}$ arrive to the system with capacity $N$ at the total rate $\lambda_k>0$, and are served at rate $\mu_k>0$ per customer. The system can be modelled as a LD-QBD $\{(X(t),\varphi(t)):t\geq 0\}$ where $X(t)$ records the total number of customers and $\varphi(t)\leq X(t)$ records the number of customers of type $1$ in the system at time $t$. The generator ${\bf Q}(\btheta)=[q(\btheta)_{(n,i)(m,j)}]$ that depends on the vector of parameters $\btheta=[\lambda_1,\lambda_2,\mu_1,\mu_2]$, is such that nonzero off-diagonals $q(\btheta)_{(n,i)(m,j)}$ are given by
	\begin{eqnarray*}
		\left\{
		\begin{array}{ll}
			\lambda_1& j = i, m=n+1, n<N;\\
			\lambda_2& j = i+1, m=n+1, n<N;\\
			(n-i)\mu_1& j = i, m=n-1, n>0;\\
			i\mu_2& j = i-1, m=n-1, n>0.
		\end{array}
		\right.
	\end{eqnarray*}
	
	Then $\frac{\partial}{\partial\btheta}{\bf Q}(\btheta)$ is given by
	\begin{eqnarray*}
		\frac{\partial q(\btheta)_{(n,i)(m,j)}}{\partial \lambda_1}&=&
		\left\{
		\begin{array}{ll}
			1& j = i, m=n+1, n<N;\\
			-1& j = i, m=n, 0\leq n<N;\\
			0& \mbox{otherwise;}	
		\end{array}
		\right.
	\end{eqnarray*}
	and
	\begin{eqnarray*}
		\frac{\partial q(\btheta)_{(n,i)(m,j)}}{\partial \lambda_2}&=&
		\left\{
		\begin{array}{ll}
			1& j = i+1, m=n+1, n<N;\\
			-1& j = i, m=n, 0\leq n<N;\\
			0& \mbox{otherwise;}	
		\end{array}
		\right.
	\end{eqnarray*}
	and
	\begin{eqnarray*}
		\frac{\partial q(\btheta)_{(n,i)(m,j)}}{\partial \mu_1}&=&
		\left\{
		\begin{array}{ll}
			(n-i)& j = i, m=n-1, n>0;\\
			-(n-i)\mu_1& j = i, m=n, 0<n\leq N;\\
			0& \mbox{otherwise;}	
		\end{array}
		\right.
	\end{eqnarray*}
	and
	\begin{eqnarray*}
		\frac{\partial q(\btheta)_{(n,i)(m,j)}}{\partial \mu_2}&=&
		\left\{
		\begin{array}{ll}
			i& j = i-1, m=n-1, n>0;\\
			-i& j = i, m=n, 0<n\leq N;\\
			0& \mbox{otherwise.}	
		\end{array}
		\right.
	\end{eqnarray*}
\end{Example}

\begin{Example}
 Suppose that the generator of a LD-QBD is a function of $\bepsilon=[\epsilon_i]_{i=1,\ldots,k}>{\bf 0}$ such that
	\begin{eqnarray*}
	{\bf Q}(\bepsilon)&=&
	{\bf Q} + \sum_{i=1}^k \epsilon_i \times \widetilde{\bf Q}_i 
	\end{eqnarray*}
	is a generator for sufficiently small $||\bepsilon||>{\bf 0}$. Then 
	\begin{eqnarray*}
		\frac{\partial {\bf Q}(\bepsilon)}{\partial \bepsilon}
		&=&
		\left[
		\begin{array}{ccc}
			\frac{\partial {\bf Q}(\bepsilon)}{\partial \epsilon_1}
			& \ldots
			& \frac{\partial {\bf Q}(\bepsilon)}{\partial \epsilon_k}
		\end{array}
		\right]
		=
		\left[
		\begin{array}{ccc}
			\widetilde{\bf Q}_1
			& \ldots
			& \widetilde{\bf Q}_k
		\end{array}
		\right].
	\end{eqnarray*}
\end{Example}

The derivatives $\frac{\partial}{\partial\btheta}$ of quantities of interest for these and other LD-QBDs can be expressed in terms of $\frac{\partial}{\partial\btheta}{\bf Q}(\btheta)$ using expressions from the matrix calculus e.g.~\cite{2007PK}, as follows. 

Let ${\bf G}^{n,n-k}={\bf G}^{n,n-k}(0)$ and $\mathbb{E}^{n,n-k}=-\frac{\partial}{\partial s}{\bf G}^{n,n-k}(s)\big|_{s=0}$ be the probability and expectation matrix, respectively. 

Let ${\bf G}^{n,n-k}(\btheta)$ and $\mathbb{E}^{n,n-k}(\btheta)$ be the notation for ${\bf G}^{n,n-k}$ and $\mathbb{E}^{n,n-k}$ when evaluated for a given $\btheta$. 

By the recursive expressions in Section~\ref{sec:Sojtimes}, we have,
\begin{eqnarray*}
	\lefteqn{
	\frac{\partial }{\partial \btheta}{\bf G}^{N,N-1}(\btheta)
	}
	\nonumber\\
	&=&
	({\bf Q}^{[N,N]}(\btheta))^{-1}
	\times
	\frac{\partial {\bf Q}^{[N,N]}(\btheta)}{\partial \btheta}
	\times
	\left(
	{\bf I}_k
	\otimes ({\bf Q}^{[N,N]}(\btheta))^{-1}
	\right)
	\\
	&&
	\times\left(
	{\bf I}_k\otimes {\bf Q}^{[N,N-1]}(\btheta)
	\right)
	\nonumber\\
	&&
	-
	({\bf Q}^{[N,N]}(\btheta))^{-1}
	\times
	\frac{\partial {\bf Q}^{[N,N-1]}(\btheta)}{\partial \btheta};
\end{eqnarray*}
for $n=N-1,\ldots,1$, we have the recursion
\begin{eqnarray*}
	\lefteqn{
	\frac{\partial }{\partial \btheta}{\bf G}^{n,n-1}(\btheta)
	}
	\\
	&=&
	-
	\frac{\partial }{\partial \btheta}
	\left({\bf Q}^{[n,n]}(\btheta)
	+{\bf Q}^{[n,n+1]}(\btheta){\bf G}^{n+1,n}(\btheta)
	\right)^{-1}
	\\
	&&
	\times\left(
	{\bf I}_k\otimes {\bf Q}^{[n,n-1]}(\btheta)
	\right)
	\\
	&&
	-
	({\bf Q}^{[n,n]}(\btheta)
	+{\bf Q}^{[n,n+1]}(\btheta){\bf G}^{n+1,n}(\btheta)
	)^{-1}
	\\
	&&
	\times
	\frac{\partial {\bf Q}^{[n,n-1]}(\btheta)}{\partial \btheta}
\end{eqnarray*}
with
\begin{eqnarray*}
\lefteqn{
	\frac{\partial }{\partial \btheta}
\left({\bf Q}^{[n,n]}(\btheta)
+{\bf Q}^{[n,n+1]}(\btheta){\bf G}^{n+1,n}(\btheta)
\right)^{-1}
}
\nonumber\\
&=&
	\left({\bf Q}^{[n,n]}(\btheta)
	+{\bf Q}^{[n,n+1]}(\btheta){\bf G}^{n+1,n}(\btheta)
	\right)^{-1}
\nonumber\\
&&
\times
\Big(
\frac{\partial ({\bf Q}^{[n,n]}(\btheta)
}{\partial \btheta}
+
\frac{\partial {\bf Q}^{[n,n+1]}(\btheta)}{ \partial \btheta}
\times\left(
{\bf I}_k\otimes {\bf G}^{n+1,n}(\btheta)
\right)
\\
&&
+
{\bf Q}^{[n,n+1]}(\btheta)
\times
\frac{\partial {\bf G}^{n+1,n}(\btheta)}{\partial \btheta}
\Big)
\nonumber\\
&&
\times
\left(
{\bf I}_k
\otimes (({\bf Q}^{[n,n]}(\btheta)
+{\bf Q}^{[n,n+1]}(\btheta){\bf G}^{n+1,n}(\btheta))^{-1}
\right)
;
\label{eq:term1}
\end{eqnarray*} 
and so for $k\geq 2$ we obtain the recursion,
\begin{eqnarray*}
	\lefteqn{
	\frac{\partial}{\partial\btheta}{\bf G}^{n,n-k}(\btheta)
}
	\\
	&=&
	\frac{\partial {\bf G}^{n,n-k+1}(\btheta)}{ \partial \btheta}
	\times\left(
	{\bf I}_k\otimes {\bf G}^{n-k+1,n-k}(\btheta)
	\right)
	\\
	&&
	+
	{\bf G}^{n,n-k+1}(\btheta)
	\times
	\frac{\partial {\bf G}^{n-k+1,n-k}(\btheta)}{\partial \btheta}.
\end{eqnarray*}
We apply similar methods to evaluate $\frac{\partial}{\partial\btheta}\mathbb{E}^{n,n-k}(\btheta)$ and the derivatives of the higher moments. 
The expressions for $\frac{\partial}{\partial\btheta}{\bf H}^{n,n+k}(\btheta)$ and related quantities follow by symmetry.

Next, to evaluate $\frac{\partial}{\partial\btheta}
{\bf f}(t;\btheta)$, we apply
\begin{eqnarray*}
\int_0^{\infty}e^{-st}
\frac{\partial}{\partial\btheta}
{\bf f}(t;\btheta)dt
&=&
\frac{\partial}{\partial\btheta}
\int_0^{\infty}e^{-st}
{\bf f}(t;\btheta)dt
=
\frac{\partial}{\partial\btheta}
\widetilde{\bf f}(s),
\end{eqnarray*}
since then the right-hand side can be computed using the results from the earlier sections, and then inverted to obtain the quantities of interest, $\frac{\partial}{\partial\btheta}
{\bf f}_n(t;\btheta)$, for all $n$.

By the recursive expressions in Section~\ref{sec:Distrt}, we have,
\begin{eqnarray*}
	\lefteqn{
	\frac{\partial \widetilde{\bf f}_{n_0}(s)}{ \partial \btheta}
=
	-\widetilde{\bf f}_{n_0}(s)
	\times
	\frac{\partial }	
	{\partial \btheta}
		\Big(
		({\bf Q}^{[n_0,n_0]}-s{\bf I})
}
	\\
	&&
		+
		\widehat {\bf R}^{(n_0-1)}(s){\bf Q}^{[n_0-1,n_0]}
		+
		\widetilde{\bf R}^{(n_0+1)}(s){\bf Q}^{[n_0+1,n_0]}
		\Big)
	\\
	&&
	\times
	\Big(
	{\bf I}_k\otimes 
	\Big(
	({\bf Q}^{[n_0,n_0]}-s{\bf I})
	\\
	&&
	+
	\widehat {\bf R}^{(n_0-1)}(s){\bf Q}^{[n_0-1,n_0]}
	+
	\widetilde{\bf R}^{(n_0+1)}(s){\bf Q}^{[n_0+1,n_0]}
	\Big)
	\Big)^{-1};
\end{eqnarray*}
and for $n\not= n_0$,
\begin{eqnarray*}
\lefteqn{
	\frac{\partial 
		\widetilde{\bf f}_{n}(s;\btheta)
	}{ \partial \btheta}
=
	\Big(
	-\bm{\alpha}
	\frac{\partial 
		{\bf W}^{n_0,n}(s;\btheta)	
	}{ \partial \btheta}
}
\\
&&	-\widetilde{\bf f}_{n}(s;\btheta)
	\times
	\frac{\partial 
	}{\partial \btheta}
	\Big(
	({\bf Q}^{[n,n]}(\btheta)-s{\bf I})
	\\
	&&
	+
	\widehat {\bf R}^{(n-1)}(s;\btheta){\bf Q}^{[n-1,n]}(\btheta)
	+
	\widetilde{\bf R}^{(n+1)}(s;\btheta){\bf Q}^{[n+1,n]}(\btheta)
	\Big)
	\\
	&&
	\times
	\Big(
	{\bf I}_k\otimes 
	\Big(
	({\bf Q}^{[n,n]}(\btheta)-s{\bf I})
\\
&&	+
	\widehat {\bf R}^{(n-1)}(s;\btheta){\bf Q}^{[n-1,n]}(\btheta)
	+
	\widetilde{\bf R}^{(n+1)}(s;\btheta){\bf Q}^{[n+1,n]}(\btheta)
	\Big)^{-1}
	\Big).
\end{eqnarray*}

Further, by the recursive expressions in Section~\ref{sec:StatDistr},
\begin{eqnarray*}
	\frac{\partial \bm{\pi}_{n}(\btheta)}{\partial\btheta}
	&=&
	\frac{\partial \bm{\pi}_{n+1}(\btheta)}{ \partial \btheta}
	\times\left(
	{\bf I}_k\otimes \widehat{\bf R}^{(n)}(\btheta)
	\right)
		\\
		&&
	+
	\bm{\pi}_{n+1}(\btheta)
	\times
	\frac{\partial \widehat{\bf R}^{(n)}(\btheta)}{\partial \btheta},
\end{eqnarray*}
where,
\begin{eqnarray*}
\lefteqn{
	\frac{\partial \bm{\pi}_N(\btheta)}{ \partial \btheta}
=
	-
	\bm{\pi}_N(\btheta)
}
\\
&&
	\times
	\frac{\partial 
	}
	{\partial \btheta}
		\left(
		\widehat {\bf R}^{(N-1)}(\btheta)
		{\bf Q}^{[N-1,N]}(\btheta)+{\bf Q}^{[N,N]}(\btheta) 
		\right)
	\\
	&&
	\times\left(
	{\bf I}_k\otimes \left(
	\widehat {\bf R}^{(N-1)}(\btheta)
	{\bf Q}^{[N-1,N]}(\btheta)+{\bf Q}^{[N,N]}(\btheta) 
	\right)^{-1}
	\right)
	.
\end{eqnarray*}

Finally, the quantities $$\frac{\partial}{\partial\btheta}\widehat {\bf R}^{(n-1)}(s;\btheta)\quad\mbox{ and }\quad \frac{\partial}{\partial\btheta}\widetilde {\bf R}^{(n+1)}(s;\btheta)$$ on the right-hand side in the above, which are required to complete the analysis, can be derived from the recursive expressions for $\widehat{\bf R}^{(n-1)}(s;\btheta)$ and $\widetilde{\bf R}^{(n-1)}(s;\btheta)$ in Section~\ref{sec:StatDistr} and Section~\ref{sec:Distrt}, respectively, by analogous methods.


\bibliographystyle{abbrv}
\bibliography{refs_sensitivity_QBDs}

\end{document}